%% file: main.tex
\documentclass[paper,11pt,nopagenum]{AAS}		% for final proceedings paper version
\PaperNumber{25-791}

\input{commands/commands}

\begin{document}

\title{Modeling and Control for Distributed Measurements of the Earth's Energy Imbalance}
\author{Rayan Mazouz,\
% \thanks{Jet Propulsion Laboratory, California Institute of Technology, Pasadena, USA.} 
% \footnotemark
Marco Quadrelli,
% \thanks{Email:{\color{blue}\href{mailto:marco.b.quadrelli@jpl.nasa.gov}{marco.b.quadrelli@jpl.nasa.gov}}}
% % \thanks{Robotics, Group Supervisor/Principal Member, Jet Propulsion Laboratory, California Institute of Technology, Pasadena, USA. Email:{\color{blue}\href{mailto:marco.b.quadrelli@jpl.nasa.gov}{marco.b.quadrelli@jpl.nasa.gov}}}
% \
Rashied Amini,
% \thanks{ Email:{\color{blue}\href{mailto:rashied.amini@jpl.nasa.gov}{rashied.amini@jpl.nasa.gov}}}
% \ 
Maria Hakuba, \\
% \thanks{Email:{\color{blue}\href{mailto:rashied.amini@jpl.nasa.gov}{rashied.amini@jpl.nasa.gov}}}
Charles Reynerson, and
% \thanks{Email:{\color{blue}\href{mailto:rashied.amini@jpl.nasa.gov}{rashied.amini@jpl.nasa.gov}}}
David Wiese}
% % \ , and 
% % \thanks{..., ..., NASA Jet Propulsion Laboratory. }

\maketitle{}

\input{Sections/0_abstract}
\footnotetext{Authors are with the Jet Propulsion Laboratory, California Institute of Technology, Pasadena, USA.}

\pagestyle{plain}

\input{Sections/1_introduction}
\input{Sections/2_dynamics}
\input{Sections/2b_attitude}

\input{Sections/2c_facet}
\input{Sections/3_control}

\input{Sections/4_experiments}
\input{Sections/5_conclusion}

\begingroup
\small
\section{Acknowledgment}
{ \noindent \textcopyright \hspace{1mm} 2025 California Institute of Technology. Government sponsorship acknowledged. This research was carried out at the Jet Propulsion Laboratory, California Institute of Technology, under a contract with the National Aeronautics and Space Administration.}

\bibliographystyle{apalike} 
\bibliography{references}
\endgroup

\end{document}

%% file: commands/commands.tex
\usepackage{bm}
\usepackage{amsmath}
\usepackage{subfigure}
\usepackage[colorlinks=false, pdfstartview=FitV, linkcolor=red, citecolor=red, urlcolor=red]{hyperref}
\usepackage[round]{natbib}
\usepackage{graphicx}
\usepackage{subcaption}
\usepackage{subcaption}

\usepackage{wrapfig} 
\usepackage{amssymb}

%% file: Sections/0_abstract.tex
\begin{abstract}
This paper presents a modeling and control framework for distributed systems in low Earth orbit, with the scientific objective of obtaining high accuracy estimates of the Earth's Energy Imbalance (EEI). This metric robustly quantifies the difference between the absorbed solar radiation, and the infrared radiation emitted into space. Formally, the EEI represents the globally and annually integrated net radiative flux at the top of the atmosphere.  
The EEI is directly correlated with physical variations in the Earth system. Obtaining accurate measurements hereof poses a major technological challenge, attributed to calibration errors of current spaceborn radiometers. This work presents a modeling and control framework for in-orbit EEI monitoring and mapping with high precision, using a distributed array of spherical spacecraft.
Perturbations and their effects on orbit and attitude are modeled, accounting for spacecraft shape and thermo-optical properties, and are subsequently used to derive optimal control for maintaining an appropriate spin rate. 
This enables each spacecraft to align closely with the orbital normal with coordinated attitudes across the formation, leading to improved spatiotemporal resolution in EEI estimation.
\end{abstract}

%% file: Sections/1_introduction.tex
\section{Introduction}
% \vspace{2.5mm}
The \emph{Earth's Energy Imbalance} (EEI) provides a valuable perspective for studying the planet's energy system. \cite{loeb2021satellite} describe how most of the EEI heats up the ocean, while the remainder heats the land or warms the atmosphere. 
Two independent studies, one based on satellite imagery and the other based on in-situ measurements, conclude that the last two decades have shown a significant increase in the rate of energy uptake, attributed to decreased reflection back into space. 
\cite{chenal2022observational} state that rigorous EEI  measurements aid in constraining the estimates of radiative sensitivity.
High-precision measurements of this parameter are hence of great value, which remains a great challenge given current Earth observation technologies.

The mission concept introduced in \cite{hakuba2023measuring,hakuba2024measuring} aims to obtain high accuracy estimates of the EEI. The proposed method is based on accelerometry that measures non-gravitational forces acting on a spacecraft in orbit. While this has been considered in the past, current state-of-the-art accelerometers allow for measurements of the EEI at an insufficient precision. The technological challenge relates directly to external disturbances, making the interpretation of radiation
flux-induced accelerations challenging. To aid the interpretation and achieve the desired accuracy, multiple spacecraft in orbit are proposed, specifically to fill the gap in diurnal sampling and attain full global coverage \citep{hakuba2024modeling}. 
The uncertainty gap is nonetheless still substantial for deriving absolute net radiative flux magnitudes. The motivation for the proposed mission naturally follows: develop a high-accuracy model of a spacecraft formation mission to obtain more comprehensive EEI estimates, serving a complementary purpose to existing methods.

A key finding by \cite{hakuba2024modeling} is that the shape and thermo-optical properties of the spacecraft establish a crucial relationship between the perturbations and the accelerations. In fact, the dependency on the cross-sectional area and incident angle inhibit a direct co-linear relationship. To accurately resolve the spatially varying incidence this work utilizes a multi-facet approach, enabling accurate modeling of the perturbative forces essential for reliable EEI estimation.
Moreover, a \emph{high LEO orbit} is preferred in which confounding aerodynamic force effects are minimized \citep{boudon1986measurement}, while remaining beneath the impact of the Van Allan belts \citep{shprits2018dynamics}. To obtain a near full sample of the globe, a near-circular sun-synchronous orbit is proposed. 
A formation of three-six spacecraft achieve a one-hour sampling resolution \citep{haar1979theoretical}. Each spacecraft must have a controlled spin rate to maintain a uniform temperature distribution throughout the mission, as well as equalize shape and optical non-uniformity. 

This paper builds upon the fundamentals established by \cite{hakuba2024modeling} and \cite{mazouz_dynamics_2021}, specifically focusing on \emph{modeling of a controlled distributed network} of spacecraft. The orbital and attitude dynamics modules accounts for perturbations, which include geopotential, atmospheric drag, solar radiation pressure, lunar and solar gravitational accelerations, and the Earth albedo. Relativistic corrections are likewise accounted for. This is an extension of the \textit{Distributed Element Beamformer Radar for Ice and Subsurface Sounding} mission \citep{haynes2021debris, bienert2021debris, mazouz_dynamics_2021, apa2022}.
The first novelty of this work lies in the use of a multi-facet attitude dynamics model that resolves force and torque contributions across discretized surface elements, enabling accurate characterization of radiation-induced perturbations essential for EEI estimation.
Furthermore, this work proposes a safety guided optimal control approach,  inspired by methods for docking \citep{zhang2022trajectory, 
laouar2024feasibility}, 
reaction wheel control \citep{mazouz2022dynamics},
and entry, descent, and landing \citep{acikmese2007convex,  mazouz2021convex}. Further details are delineated in the remainder of the paper.

%% file: Sections/2_dynamics.tex
\section{Dynamics}
\label{sec:dynamics}

This section first outlines the appropriate reference frames, after which the governing \emph{orbital} and \emph{attitude} equations of motion are delineated -- in light of a single  spacecraft in LEO. 
There are six predominant external accelerations that are considered in these models, which are perturbations due to the 
geopotential, atmospheric drag, solar radiation pressure, lunar and solar gravitational accelerations, and the Earth albedo. 
High precision modeling also requires accounting for Earth tides accelerations and relativistic corrections \citep{montenbruck2014models}, which are likewise incorporated. 
Finally, the multi-facet attitude dynamics approach for a near-spherical spacecraft is delineated to capture differential exposure to radiation with high accuracy.

\subsection{Reference Frames}
\begin{figure}[t!]
\centering
\includegraphics[width=.40\textwidth]{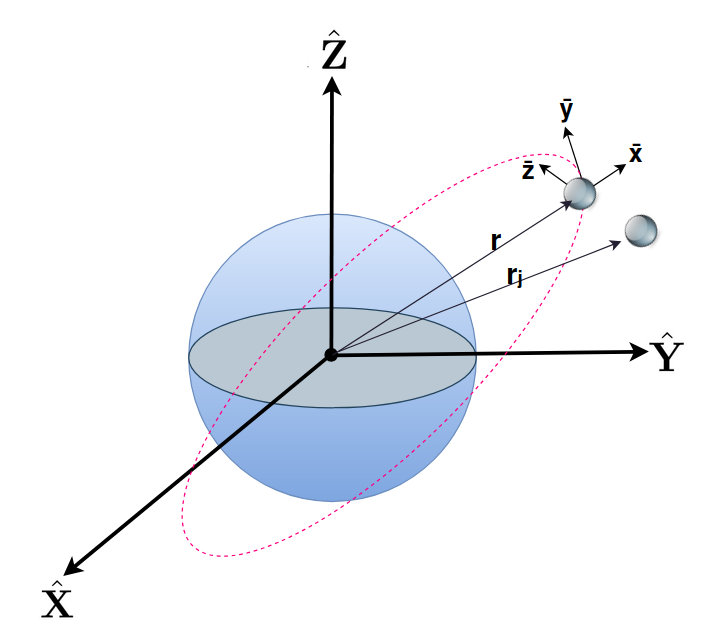}
\caption{A representation of the inertial $\mathcal{I}$- and the LVLH $\mathcal{L}$-frame. Vector $\mathbf{r}$ defines the radius of the leader spacecraft, while $\mathbf{r_j}$ specifies the radius of the $j$-th sphere, both in the inertial frame.}
\label{fig:ecilvlh}
\end{figure}

The necessary reference frames for the dynamics formulations are depicted Fig. \ref{fig:ecilvlh}, and are summarized as follows
\begin{itemize}
\item[$\mathcal{I}$]  The \emph{Earth-Centered Inertial} (ECI) frame has axis $\hat{\mathbf{X}}$ pointing towards the Vernal Equinox and $\hat{\mathbf{Z}}$ towards the North Pole, with $\hat{\mathbf{Y}}$ completing the right-handed frame.
\item[$\mathcal{L}$] The \emph{Local Vertical Local Horizontal} (LVLH) frame is attached and aligned to the leader, and has axis $\bar{\mathbf{x}}$ pointing along the local vertical (radial), $\bar{\mathbf{y}}$ towards the direction of flight (along-track), and $\bar{\mathbf{z}}$ in the orbit normal direction (across-track). 
\item[$\mathcal{B}$] The \emph{body frame} is fixed to the vehicle. The angular rate of the vehicle, expressed as the vector {\boldmath{$\omega$}}, is resorted to for attitude quantification.
\end{itemize}

\subsection{Orbital Dynamics}
 A non-linear and non-circular orbital dynamics model is considered, governed by
\begin{equation}
\label{governing}
    {\bf{\ddot{r}}} + \frac{\mu}{{r^3}} {\bf{r}} = \frac{\partial R}{\partial \mathbf{r}},
\end{equation}
in which $\frac{\partial R}{\partial \mathbf{r}} $ defines some perturbing potential $R$, and $\mathbf{r}$ the position vector in the Earth-Centered Inertial (ECI) frame. The perturbations in this reference frame are outlined subsequently. 

\vspace{-5mm}

\subsubsection{Gravity Acceleration} 

The acceleration due gravity is given by the potential at a given point, 
defined by the evaluation of Legendre polynomials $P_{nm}$
\citep{montenbruck2014models,cunningham1970computation}. 
The dependency is governed by the geocentric latitude $\phi$ and 
longitude $\lambda$ of location $\mathbf{r}$.
As outlined in \cite{milani1987}, the acceleration due to the Earth’s gravity potential is as follows
\begin{align}
    \ddot{\mathbf{r}}_{\text{grav}} = \nabla \frac{\mu}{r} 
    \sum_{n=0}^{\infty}
    \sum_{m=0}^{\infty}
    \frac{Re^n}{r^n}
    \bar{P}_{nm}
    \sin(\phi)
    (
        \bar{C}_{nm} \cos(m \lambda)  +
        \bar{S}_{nm} \sin(m \lambda) 
    ).
\end{align}
Note that the coefficients $\bar{C}_{nm}$ and $\bar{S}_{nm}$ describe the normalized 
dependence on the Earth’s internal mass distribution. 
For the computation of these coefficients, along with normalized associated 
Legendre function $\bar{P}_{nm}$, 
the reader is referred to \cite{montenbruck2014models}. Term expansion based on zonal and tesseral-sectorial geopotential is described in \cite{mazouz_dynamics_2021}.
\vspace{-5mm}

\subsubsection{Atmospheric Drag}

Atmospheric drag represents one of the largest
non-gravitational perturbations acting on the spacecraft in LEO \citep{montenbruck2014models,mazouz_dynamics_2021}.
Let $m$ and $A$ denote the mass and surface area of the spacecraft, respectively, such that
\begin{align}
\label{eq:drag}
    \ddot{\mathbf{r}}_{\text{drag}} = -\frac{1}{2} C_D \frac{A}{m} \rho ({\bf{v}} -{\bf{v}}_{\text{atm}}) \left \| ({\bf{v}} -{\bf{v}}_{\text{atm}}) \right \|,
\end{align}
defines the acceleration due to drag.
Further, $C_{D}$ is the drag coefficient and $\rho$ the local atmospheric density. 
Atmospheric properties are approximated using the Jacchia-Robert drag model, based on an analytic model for upper atmospheres \citep{roberts_analytic_1971}. 
\vspace{-5mm}

\subsubsection{Solar Radiation Pressure}

The solar flux determines the magnitude of the solar radiation pressure. 
Accounting for the inclined angles of the incoming radiation, the acceleration is defined by
\begin{align}
\label{eq:srp}
    \mathbf{\ddot{r}}_{\text{srp}} = - \eta P_{\odot} \frac{1 AU^2}{r_{\odot}^2}
    \frac{A}{m}
    \cos \theta
    \left [
    (1-\epsilon)\mathbf{e}_{\odot}
    +
    2\epsilon \cos \theta \mathbf{n}
    \right ].
\end{align}
Note that $\epsilon$ denotes the reflectivity, and $\mathbf{e}_{\odot}$, 
$\mathbf{n}$ are unit vectors, such that $\cos \theta = 
\mathbf{n}^{T} \mathbf{e}_{\odot}$. To enhance the accuracy of the solar 
radiation pressure acceleration, an eclipse condition and shadow functions are incorporated,
% The eclipse condition assumes a conical shadow model, where oblateness is neglected.
established by parameter $\eta$, such that 
\begin{align*}
    \eta =
    \begin{cases}
      0, & \text{if spacecraft is in umbra,}\  \\
      1, & \text{if spacecraft is in light conditions,}\ \\
      (0,1) & \text{if spacecraft is in penumbra}.
    \end{cases}
\end{align*}
\vspace{-10mm}

\subsubsection{Third Body Accelerations}

The lunar and solar perturbing accelerations are considered. Given that the bodies 
are much further away from Earth than the spacecraft, the accelerations follow from
\begin{align} 
    \ddot{\mathbf{r}}_{\text{body}} = -\frac{\mu}{s^3} r,
\end{align}
where $s$ is the geocentric position of the perturbing body.
% , and $\mathbf{r}$ the geocentric
% coordinates of the spacecraft. 
Chebyshev approximations \citep{hernandez2001chebyshev} are utilized to approximate a solution of the nonlinear operator. 
% equation in a Banach space.
% Precise advanced analytical methods are outlined in \citep{flandern1979}. 
% A series of solar system ephemerides in the form of these approximations are available
% \citep{seidelmann1992}. This standard of high-precision approximations are incorporated into the model.
\vspace{-5mm}

\subsubsection{Earth Tides}

Due to the gravitational attraction of the Sun and the Moon exerted on the Earth, a time-varying
deformation is encountered within the internal structure of the planet.
This acceleration is categorized into two types: small periodic deformations of the solid 
body, and ocean tides responses.
% This implies that the gravity field of Earth is not static. 
This causes a shift in the earlier outlined spherical harmonics coefficients $C_{nm}, S_{nm}$. The 
expressions are derived in \cite{sanchez1976}, and adapted in this study.
\vspace{-5mm}

\subsubsection{Relativistic Effects}

Due to the mass and angular momentum of the Earth, potential fields are generated that lead to a curvature of space-time. In the vicinity of the Earth, the geodesic equation is approximated 
using first-order expansion \citep{mccarthy1996iers}, such that
\begin{align}
    \mathbf{\ddot{r}}_{\text{rel}} =  
    +
    \frac{\mu}{r^2}
    \left (
    (\frac{4 \mu}{c^2 r} - \frac{v^2}{c^2})\mathbf{e}_r
    + \frac{4v^2}{c^2}(\mathbf{e}_r \cdot \mathbf{e}_v)\mathbf{e}_v 
    \right),
\end{align}
where $c$ and $v$ denote the speed of light and the spacecraft, respectively. Further, $\mathbf{e}_r, \mathbf{e}_v$ denote the unit position and velocity vectors. This acceleration discards gravito-magnetic interactions. 
\vspace{-5mm}

\begin{figure}[t!]
\centering
\includegraphics[width=.65\textwidth]{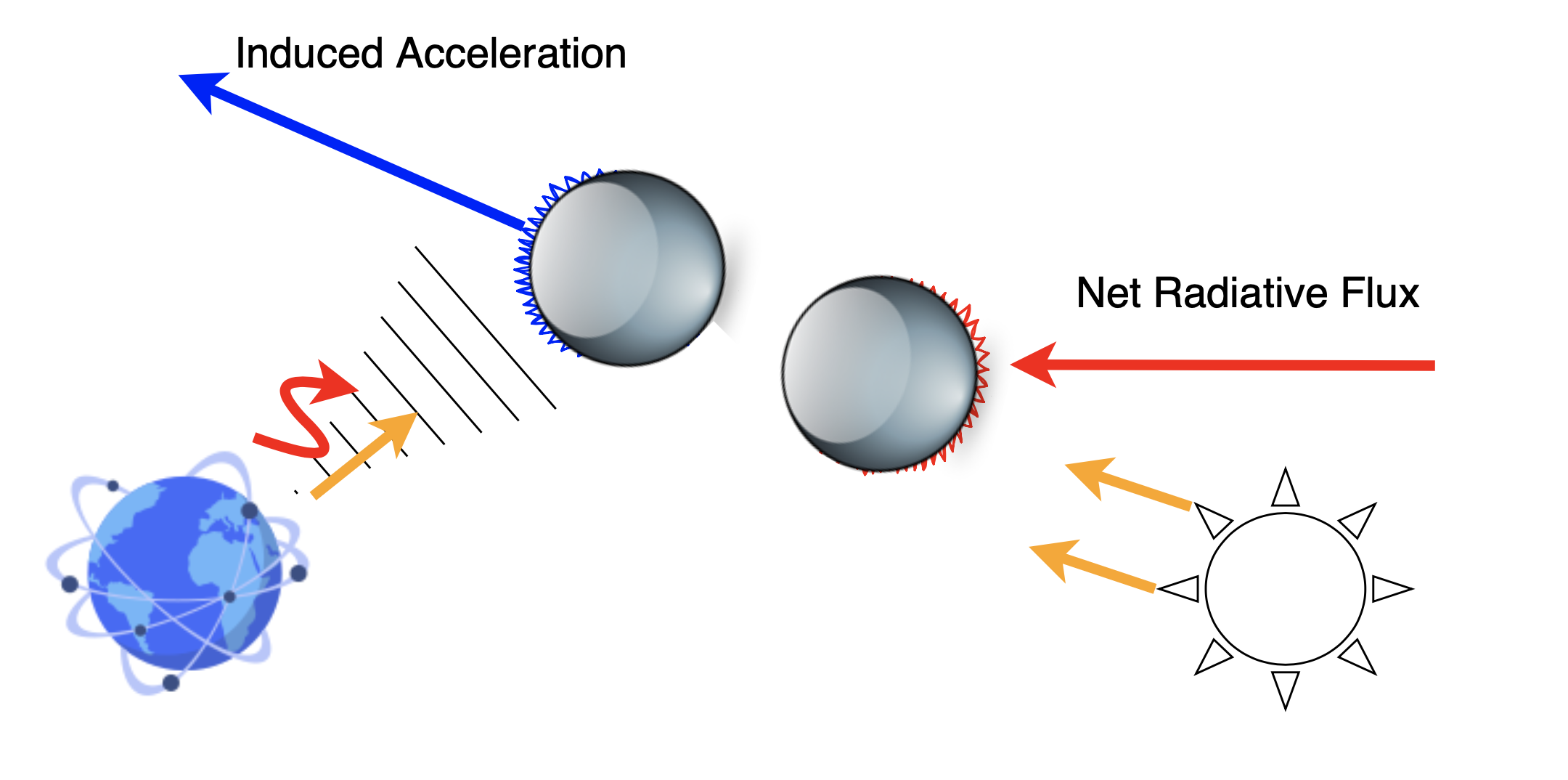}
\caption{Radiation pressure acting on spherical spacecraft. The acceleration is function of the effective cross-sections and the impinging irradiance.}
\label{fig:earth_albdeo}
\end{figure}

\subsubsection{Earth Albedo}

This acceleration is central to the application at hand, of which the conceptual idea depicted in Fig. \ref{fig:earth_albdeo}.
For a perfectly spherical spacecraft with constant geometric cross section, and well-known optical coating property coefficient, acceleration and irradiance are co-linear and proportional \citep{hakuba2024measuring}
The total acceleration on a spherical satellite of constant cross sectional area and invariant reflective and emissive properties can be approximated as
\begin{equation}
\label{eq:total_albedo}
   \mathbf{\ddot{r}}_{\text{alb}} = 
    \sum_{j=1}^{N}
    \text{d}\phantom{.} \overrightarrow{a}_j,
\end{equation}
where
\begin{equation}
\label{eq:ind_albedo}
  \text{d}\phantom{.} \overrightarrow{a}_j
  =
  K
  \left [
  (\nu a E_s \cos \theta_{s}
  + e M_{B})
  \frac{A_c}{m c \pi r^2}
  \cos \alpha
  ~
  dA~
  \hat{r}
    \right ]_{j}.
\end{equation}
Parameter $K = (1 + \eta_{E})$ denotes the Lochry radiation augmentation factor, modeled as function of the solar radiation satellite reflectivity $\eta_{E}$. Indicator function $\nu$ specifies whether a given Earth element is in light ($\nu = 1$) or darkness ($\nu = 0$). Next, $\Phi_{\text{IN}} = E_{S} \cos \theta_{s}$ defines the total amount of solar flux incident on a (planar) earth element. Here, $E_{s}$ is the solar irradiance (in  Wm$^{-2}$ ) for an isotropic Sun, and $\theta_{s}$ the solar zenith angle. The exitance of $dA$ is defined by $M_b$, and by assuming the Earth is an ideal black body emitting incident solar radiation isotropically, it follows that $M_b = \frac{E_s}{4}.$
Term $\frac{A_c}{m}$ defines the area-to-mass ratios, while terms $a$ and $e$ define the  albedo and emissivity models, respectively, defined as follows
\begin{align*}
    a = a_o + a_1 P_1(\sin \phi) + a_2 P_2(\sin \phi),  \\
    e = e_o + e_1 P_1(\sin \phi) + e_2 P_2(\sin \phi), 
\end{align*}
where $\phi$ defines the equatorial latitude.
As in \cite{knocke1988earth} using the zeroeth, first and second degree zonal harmonics suffices. Finally, the projected attenuated area is 
the element area $dA$ multiplied by the cosine of the view angle $\alpha$
\begin{align*}
    A' = \frac{dA_j \cos \alpha_j}{\pi r^2} \hat{r},
\end{align*}
where each area $dA_j$ corresponds to the surface 
portion of a sphere with a mean Earth radius.

%% file: Sections/2b_attitude.tex
\subsection{Attitude Dynamics}

Next, the attitude dynamics formulation is outlined, which is specifically important for the analysis of external perturbation forces, requiring precise modeling to extract the contribution of the EEI. By the transport theorem, the governing equation for the attitude is
\begin{equation}
[J_B] \dot{\omega}_{B/I} + \omega_{B/I} \times [J_B] \omega_{B/I} =
\sigma + \tau.
\end{equation}
In this equation, \(\omega_{B/I}\) represents the spacecraft body rate vector relative to an inertial frame, while \([J_B]\) denotes the spacecraft's inertia tensor, which characterizes its resistance to changes in rotational motion. The term \(\sigma\) accounts for the sum of external perturbations, while \(\tau\) represents the control torque applied to achieve desired orientation and stability. For a quaternion \(\mathbf{q} = [q_1, q_2, q_3, q_4]^T\) parameterization, where \(q_4\) is scalar and \([q_1, q_2, q_3]^T\) a vector, the quaternion kinematic equation is
\begin{equation}
\label{eq:kinematic}
\dot{\mathbf{q}} = \frac{1}{2} \mathbf{q} \otimes \begin{bmatrix} \omega_{B/I} \\ 0 \end{bmatrix}.
\end{equation}
In the remainder of this section the terms contributing to $\sigma$ are outlined. 
\vspace{-5mm}

\subsubsection{Gravity Gradient Torque}

The gravity gradient torque arises due to the interaction between Earth's gravitational field and the spacecraft's mass distribution
\[
\sigma_{\text{grav}} = 3 \frac{\mu}{r^3} \mathbf{r} \times \left[ [J_B] \mathbf{r} \right],
\]
where \(\mu\) is Earth's gravitational parameter. If the spacecraft's inertia tensor is isotropic (i.e., a sphere), \([J_B]\mathbf{r}\) is aligned with \(\mathbf{r}\), resulting in \(\sigma_{\text{grav}} = 0\). However, given that the spacecraft may have slight misalignments or asymmetries, this torque is not discarded.
\vspace{-5mm}

\subsubsection{Drag Torque}

Atmospheric drag torque arises from the interaction between the residual atmospheric particles and the spacecraft's surface, creating a force offset from the center of mass (CM). Recall the acceleration of atmospheric drag given in Eq.~\eqref{eq:drag}. Let $\mathbf{r}_{cp}$ be the non-zero position vector of the center of pressure (CP), relative to the center of mass of the body, such that
\[
\sigma_{\text{drag}} = \mathbf{r}_{\text{CP}} \times (m \cdot {\ddot{\mathbf{r}}_{\text{drag}}}).
\]
This torque depends on the spacecraft's shape, orientation, and atmospheric conditions. For a spherical spacecraft, the CP typically coincides with the CM, if the sphere is symmetric and the drag is uniformly distributed. However, slight asymmetries create an offset, which are accounted for. 
\vspace{-5mm}

\subsubsection{Solar Pressure Torque}

Solar radiation pressure torque arises from the force exerted by solar photons on the spacecraft's surface. The corresponding torque is defined as
\[
\sigma_{\text{srp}} = \mathbf{r}_{\text{CP}} \times (m \cdot \ddot{\mathbf{r}}_{\text{srp}}),
\]
where \(\ddot{\mathbf{r}}_{\text{srp}}\) is defined in Eq.~\eqref{eq:srp}. Not all facets are equally affected by the acceleration due to solar radiation pressure, resulting in a non-uniform distribution of forces that generates a net torque. 
\vspace{-5mm}

\subsubsection{Earth Albedo Torque}

The final predominant torque considered for modeling the attitude dynamics is the torque due to the Earth's albedo. The torque arises from the reflection of solar radiation by the Earth's surface. The total acceleration caused by this effect is given by Eq.~\eqref{eq:total_albedo}, where each surface element contributes to the overall force, and is described by Eq.~\eqref{eq:ind_albedo}.  A graphical depiction has been provided in Fig. \ref{fig:earth_albdeo}. 
Specifically, each element \(j\) experiences a force that depends on factors such as the spacecraft's reflectivity, the solar irradiance, and the angle of incidence. The attitude torque due to Earth albedo is computed as the cross product of the CP vector and the acceleration is derived from the individual contributions of each surface element
\[
\sigma_{\text{alb}} = \mathbf{r}_{\text{CP}} \times (m \cdot \ddot{\mathbf{r}}_{\text{alb}}),
\]
where \(\ddot{\mathbf{r}}_{\text{alb}}\) is given by Eq.~\eqref{eq:total_albedo}. This torque plays a crucial role in the analysis of the EEI. 

\begin{wrapfigure}[13]{r}{0.5\textwidth}
  \centering
  \vspace{-10mm}
  \includegraphics[width=.35\textwidth]{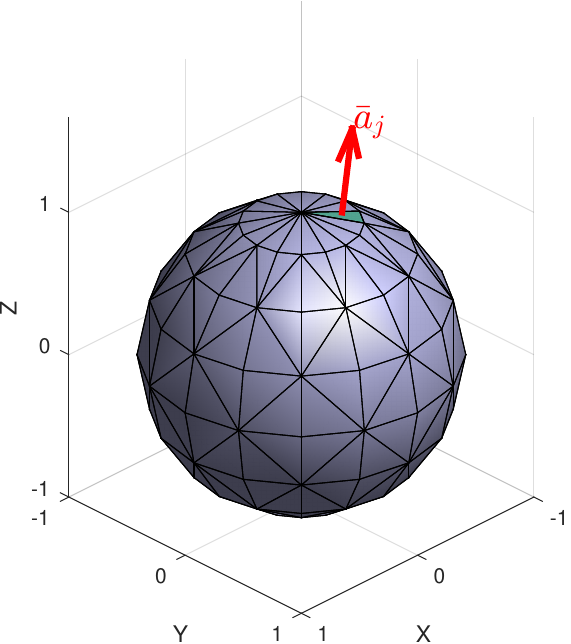}
  \caption{Spherical spacecraft object built from polyhedral facets; $\bar{a}_j$ indicates acceleration.}
  \label{fig:facet_plot}
\end{wrapfigure}

%% file: Sections/2c_facet.tex
\subsection{Multi-Facet Attitude Dynamics}

In order to accurately extract the Earth’s Energy Imbalance (EEI) from the orbital distribution, for each spherical spacecraft a multi-faceted approach is taken to compute the relevant perturbations with high fidelity.
By subdividing the spacecraft into polyhedral facets, as depicted in Fig.~\ref{fig:facet_plot}, the differential exposure of each facet to radiation can be captured accurately. Solar radiation pressure varies based on the surface area presented to the Sun, while Earth albedo radiation is reflected from Earth's surface and hits the spacecraft depending on its orbital position and attitude. The faceted approach enables a more accurate computation of the forces and torques acting on each facet, ensuring that both effects are accounted for based on the spacecraft's continuously changing orientation. 
\begin{minipage}{\textwidth}~\end{minipage}
\indent Let the spacecraft be modeled as a collection of \( N \) flat facets, each with surface area \( A_j \), unit surface normal vector \( \mathbf{n}_j \in \mathbb{R}^3 \), position vector \( \mathbf{r}_j \in \mathbb{R}^3 \) from the spacecraft's center of mass, and material reflectivity \( \rho_j \in [0,1] \). Further, let \( \mathbf{r}_\odot \in \mathbb{R}^3 \) denote the unit vector pointing from the spacecraft to the Sun, and \( \mathbf{r}_\oplus \in \mathbb{R}^3 \) the effective direction of Earth albedo. The total radiation force on facet \( j \), which results in a corresponding acceleration \( \mathbf{a}_j = \mathbf{F}_j^{\text{rad}} / m \), is modeled as
\[
\mathbf{F}_j^{\text{rad}} = \mathbf{F}_j^{\odot} + \mathbf{F}_j^{\oplus},
\]
where
\[
\mathbf{F}_j^{\odot} = - P_{\odot}  A_j \left[ (1 - \rho_j) \cos\theta_j \, \mathbf{r}_\odot + 2 \rho_j \cos\theta_j \, \mathbf{n}_j \right] \cdot \mathbb{I}_{\{\cos\theta_j > 0\}},
\]
\[
\mathbf{F}_j^{\oplus} = - P_{\oplus, j} A_j \left[ (1 - \rho_j) \cos\phi_j \, \mathbf{r}_\oplus + 2 \rho_j \cos\phi_j \, \mathbf{n}_j \right] \cdot \mathbb{I}_{\{\cos\phi_j > 0\}}.
\]
Here, \( P_{\odot}  \) denotes the solar radiation pressure at 1 AU, and \( P_{\oplus, j} \) represents the Earth albedo pressure incident on facet \( j \). The angle of incidence with respect to the Sun is given by \( \theta_j = \arccos(\mathbf{n}_j^\top \mathbf{s}_\odot) \), while the angle of incidence with respect to Earth albedo is \( \phi_j = \arccos(\mathbf{n}_j^\top \mathbf{s}_\oplus) \). \\
\begin{minipage}{\textwidth}~\end{minipage}
\indent The function \( \mathbb{I}_{\{\cdot\}} \) is the indicator function, used to ensure that only illuminated facets (those with positive cosine of the incidence angle) contribute to the net force and torque. 
The total radiation-induced torque is then given by
\[
\sigma_{\text{rad}} = \sum_{j=1}^{N} \mathbf{r}_j \times \left( \mathbf{F}_j^{\odot} + \mathbf{F}_j^{\oplus} \right).
\]
Through this  approach, modeling of the spacecraft's dynamics becomes more sophisticated. By systematically evaluating the accelerations of each facet, it is also possible to develop more effective control algorithms that account for the nuances of the spacecraft's interaction with its environment. The eventual goal is to be able to analyze these accelerations to be able to model the influence of the EEI with high precision, for which the faceted approach constitutes a more meticulous approach.

\clearpage

%% file: Sections/3_control.tex
\section{Optimal Control}
\label{sec:control}

This section introduces the design of optimal control laws for 
the distribution of spacecraft, collectively providing high-accuracy EEI measurements. The primary objectives are twofold: first, to synchronize the spin rate of each spacecraft at a minimum rate, ensuring a constant temperature; second, to maintain the spacecraft's proper orientation throughout its mission lifetime, relative to Earth and other spacecraft in the distribution. The performance of key attitude control actuators, namely magnetic torquers and reaction wheels, is analyzed in the context of these dual tasks. Through this strategy, the goal is to be able to derive the EEI in a post-processing procedure, with high accuracy. Next, the convex optimization formulation is outlined, followed by the actuator control laws.
\vspace{-3mm}

\subsection{Convex Optimization}
Let the state vector for each follower spacecraft $i = \{ 1,2, \hdots, n \}$ as $\mathbf{x}_{i}(t) \in {\mathbb{R}}^{13}$ and control vector $\mathbf{u}_{i}(t) \in {\mathbb{R}}^{6}$, such that
\begin{align}
\begin{split}
    & {\mathbf{x}}_{i}(t) = \begin{bmatrix}
    {\mathbf{r}}^{T} & {\mathbf{v}}^{T} & {\mathbf{q}}^{T} &
    {\pmb{\omega}}^{T}
    \end{bmatrix}_{i}^{T}, \qquad {\mathbf{u}}_{i}(t) = \begin{bmatrix}{u_{x}} & {u_{y}} & {u_{z}}      
    \end{bmatrix}_{i}^{T}.
     \end{split}
\end{align}
The decision vector is defined as ${\bf z}_{i}(t) = [{\mathbf{x}_{i}}(t)^{T} \quad {\mathbf{u}_{i}}(t)^{T}]^{T}$. In the formulation a leader $L$ is considered, and the $i^{\text{th}}$ follower $\mathbf{x}_i$. The following optimizations seeks to minimize a cost function $J({\bf z_{i}}(t))$,  during time-of-flight $t_{f}$, subject to
\begin{align}
\begin{split}
{\bf{\dot{x}}}_{i}(t) & = f({{\bf{x}}}, {\bf{\Gamma}}) + {\bf{B}} {\bf{u}}_{i}(t), \\
\left \| {\bf{u}}_{i}(t)\right \| & \leq u_{\text{max}}, \\
{\bf{{x}}}_{i}(t_{0}) & = {\bf{{x}}}^{0}, \\
{\bf{{x}}_{i}}(t_{f}) & = {\bf{{x}}}^{f}, \\
\left \| ({\bf{x}_{i}}(t) - {\bf{x}}_{L}(t))\right \|_{2} & \geq R_{\text{col}}.
\end{split}
\end{align}
The state \( {\bf x}(t)_{i} \) evolves according to the nonlinear dynamics \( {\bf \dot{x}}(t) = f({\bf x}, \boldsymbol{\Gamma}) + {\bf B} {\bf u}(t) \), where \( \boldsymbol{\Gamma} \) represent the orbital parameters, \( {\bf B} \) is an input matrix, and \( {\bf u}(t) \) is the control input constrained by a maximum  \( u_{\text{max}} \). The system starts at \( {\bf x}^0 \) at time \( t_0 \) and must reach a final state \( {\bf x}^f \) at time \( t_f \). Additionally, a collision avoidance constraint enforces that the distance between two entities, with positions \( {\bf x}_{i}(t) \) and \( {\bf x}_L(t) \), remains at least \( R_{\text{col}} \) apart throughout the trajectory.
\vspace{-3mm}

\subsection{Feedback Actuation}
For spacecraft actuation, two types of actuators are considered: magneto torquers and reaction wheels.
The \emph{magnetic torque} $\mathbf{T}$ exerted by magnetic torquers in the body frame is expressed as the cross product ($\mathbf{T} = \mathbf{m} \times \mathbf{B},
$) of the magnetic moment $\mathbf{m}$ of the spacecraft and the Earth's magnetic field $\mathbf{B}$. 
The magnetic moment $\mathbf{m}$ is a property of the magnetic torquers, which can be controlled to generate the desired magnetic torque. The Earth's magnetic field $\mathbf{B}$ in the body frame is typically obtained by transforming the magnetic field vector from the inertial frame to the body frame using the spacecraft's attitude information. In spherical coordinates, the Earth's magnetic field can be described with radial (\(B_r\)) and latitudinal (\(B_\theta\)) components  
\citep{gatherer2019magnetorquer}
\begin{align*}
B_r &= -2B_0 \left( \frac{R_E}{r} \right)^3 \cos \theta, \qquad
B_\theta = -B_0 \left( \frac{R_E}{r} \right)^3 \sin \theta, \qquad
|\mathbf{B}| = B_0 \left( \frac{R_E}{r} \right)^3 \sqrt{1 + 3\cos^2 \theta },
\end{align*}
where \(\theta\) is the colatitude measured from the north magnetic pole (or geomagnetic pole), and \(B_0\) the magnetic field strength at the Earth's surface.

The magnetic moment $\mathbf{m}$ is determined using a proportional-derivative (PD) control law to counteract rotational motion and maintain the desired attitude. The control law is expressed as $\mathbf{m} = -K_p \mathbf{q}_{\text{error}} - K_d (\boldsymbol{\omega}_{\text{B}} - \boldsymbol{\omega}_{\text{B,d}})$, where $K_p$ and $K_d$ are the proportional and derivative gain matrices, respectively, $\mathbf{q}_{\text{error}}$ represents the attitude error quaternion, $\boldsymbol{\omega}_{\text{B}}$ is the current angular velocity vector in the body frame, and $\boldsymbol{\omega}_{\text{B,d}}$ is the desired angular velocity vector in the body frame. This control law ensures stability and accurate pointing by compensating for attitude errors and disturbances during spacecraft operations.

Next, \emph{reaction wheels} are a popular approach towards spacecraft attitude control, by generating control torques to adjust orientation through the conservation of angular momentum. These devices work in conjunction with control strategies to achieve desired pitch, roll, and yaw orientations while compensating for external disturbances. The angular acceleration of the wheels is modeled as
$\pmb{\dot{\bar{\omega}}}_{\text{wheel}} = -\pmb{\omega} \times \pmb{\bar{\omega}}_{\text{wheel}} - \frac{1}{J_{W}} \pmb{\tau}_{\text{control}}$, where $J_{W}$ is the moment of inertia about the spin axis.

Likewise, a PD controller is employed to ensure feedback in the control loop to minimize errors. The control design is divided into pitch axis control and coupled roll-yaw axis control. For details on the latter, the reader is referred to \cite{mazouz2022dynamics}. In what follows, the experimental results are outlined, for both actuation cases.

%% file: Sections/4_experiments.tex
\section{Experiments}
\label{sec:simulate}

In this section, the simulation results are presented. First, the control effectiveness of the magnetotorquers and reaction wheels for attitude stabilization and maneuvering are compared. 
Both systems are analyzed under identical mission conditions and subject to the same external disturbances. 
For mission details (e.g., orbital elements), the reader is referred to \cite{hakuba2024measuring}.
Each multi-facet spacecraft is considered to have a total mass of $m=$ 50 kg, and a radius of $1$ m. 
\vspace{-5mm}

\subsection{Optimal Control}

Two scenarios are considered to evaluate the performance of the control
\begin{itemize}
\item[] \textit{Scenario I:} Pointing stabilization under constant external disturbances, maintaining a fixed nadir-pointing attitude.
\item[] \textit{Scenario II:} Maneuvering between two specified attitudes within a given time frame while countering dynamic external disturbances.
\end{itemize}

\begin{figure*}[t!]
\centering
    \begin{minipage}{0.49\textwidth}
        \centering
        \includegraphics[width=1.00\linewidth]{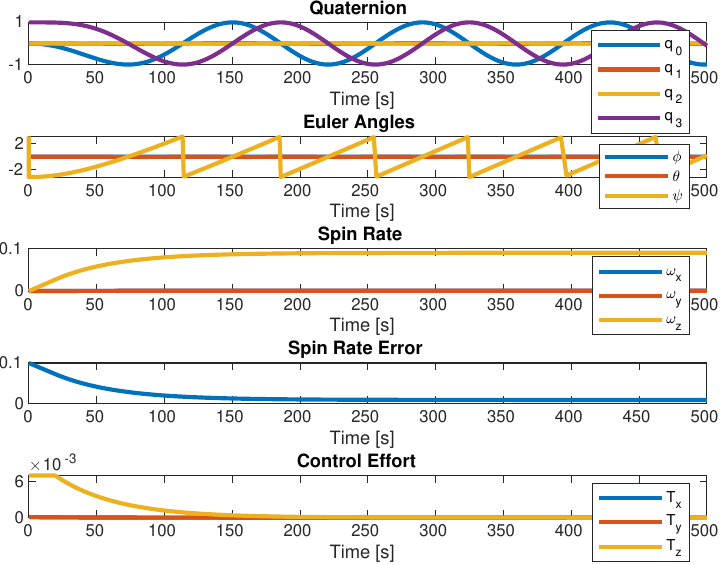}
    \end{minipage} 
    \begin{minipage}{0.50\textwidth}
        \centering
        \includegraphics[width=1.00\linewidth]{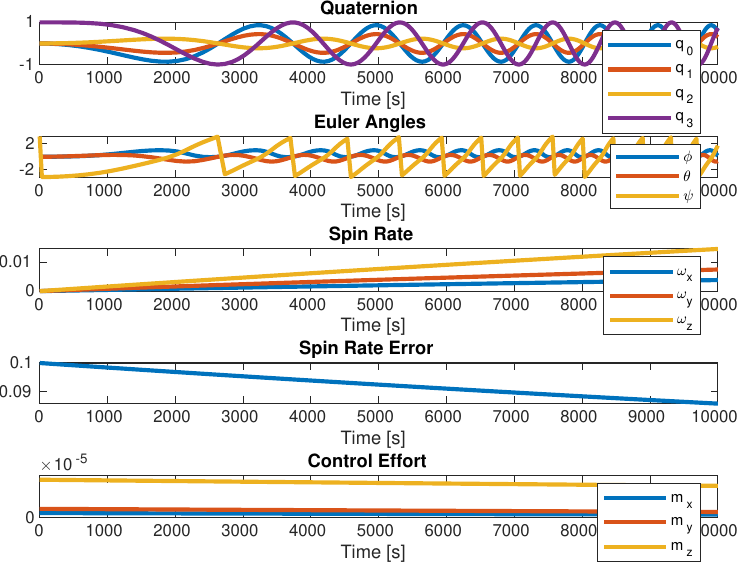}
    \end{minipage}
    \caption{Control analysis: (left) based on reaction wheels, and (right) magneto torquers based. }
\label{fig:control}
\end{figure*}

For the magneto-torquers, the magnetic field model was utilized to calculate the available control moments, while for reaction wheels, saturation and momentum dumping were considered. The trade-offs between precision and efficiency are highlighted in Fig. \ref{fig:control}, where the quaternions, Euler angles, LVLH deviation, attitude error and control effort are plotted. The results indicate that magneto-torquers are well-suited for missions with moderate pointing requirements and limited power budgets, whereas reaction wheels are preferable for high-precision applications.

The results demonstrate distinct differences in the performance of the two actuation mechanisms. Reaction wheels excel in precise attitude stabilization and maneuvering, achieving an RMS attitude error of less than $0.01^\circ$ in both scenarios. 

However, this precision comes at the cost of higher power consumption and the necessity for periodic momentum dumping.
Magnetotorquers, on the other hand, provide a simpler and more energy-efficient solution, but their performance is inherently limited by the Earth’s magnetic field strength and orientation. In Scenario I, magnetotorquers achieved an RMS attitude error of approximately $0.1^\circ$, while in Scenario II, the settling time was noticeably longer compared to reaction wheels. In fact, the magnetotorquers for all considered time-horizons are not able to obtain the desired performance.

\begin{figure*}[b!]
\centering
    \begin{minipage}{0.49 \textwidth}
        \centering
        \includegraphics[width=1.00\linewidth]{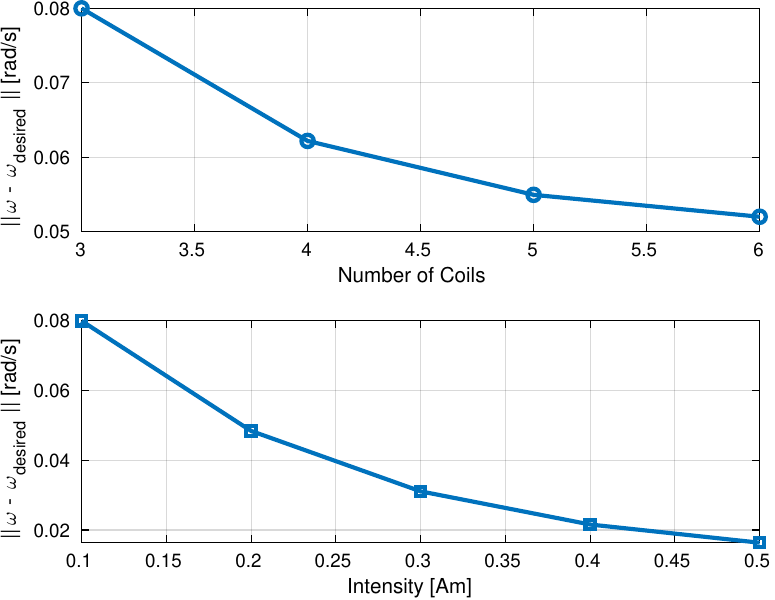}
    \end{minipage}
        \hspace{1mm}
        \begin{minipage}{0.49\textwidth}
        \centering
        \includegraphics[width=1.00\linewidth]{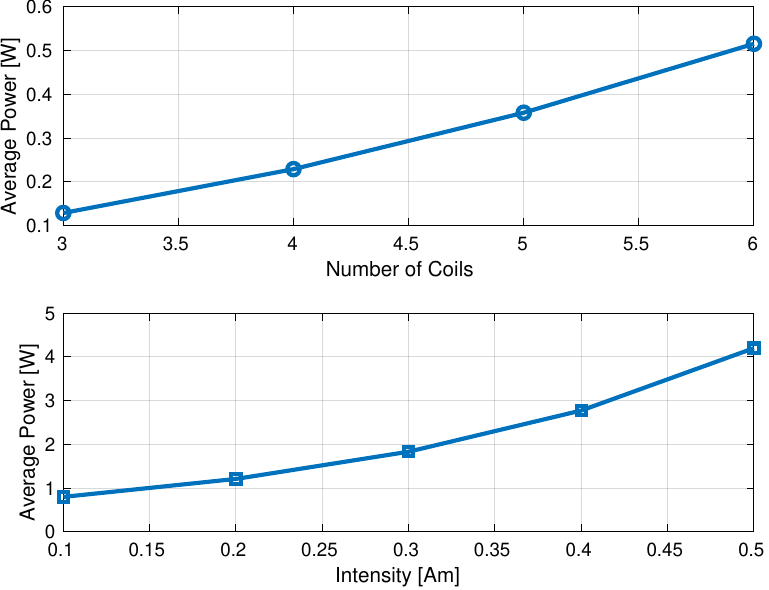}
    \end{minipage} 
       \caption{Spin rate error (left) and average power consumption (right) versus number of coils and magnetic intensity.}
\label{fig:power}
\end{figure*}

Next, it is examined how key spacecraft actuator parameters—specifically the number of magnetic coils and their respective magnetic intensity —affect the spacecraft’s ability to achieve a desired spin rate. The plots in Fig. \ref{fig:power} illustrate that increasing intensity reduces control error more effectively than simply adding coils, highlighting a trade-off between control precision and power consumption. While higher intensity improves accuracy, it comes at the expense of higher power consumption, as also depicted in the same figure.
This insight helps guide the design of magnetic torque control systems to balance efficiency and performance.

\subsection{Earth's Energy Imbalance}

\begin{figure}[t!]
\centering
\includegraphics[width=.55\textwidth]{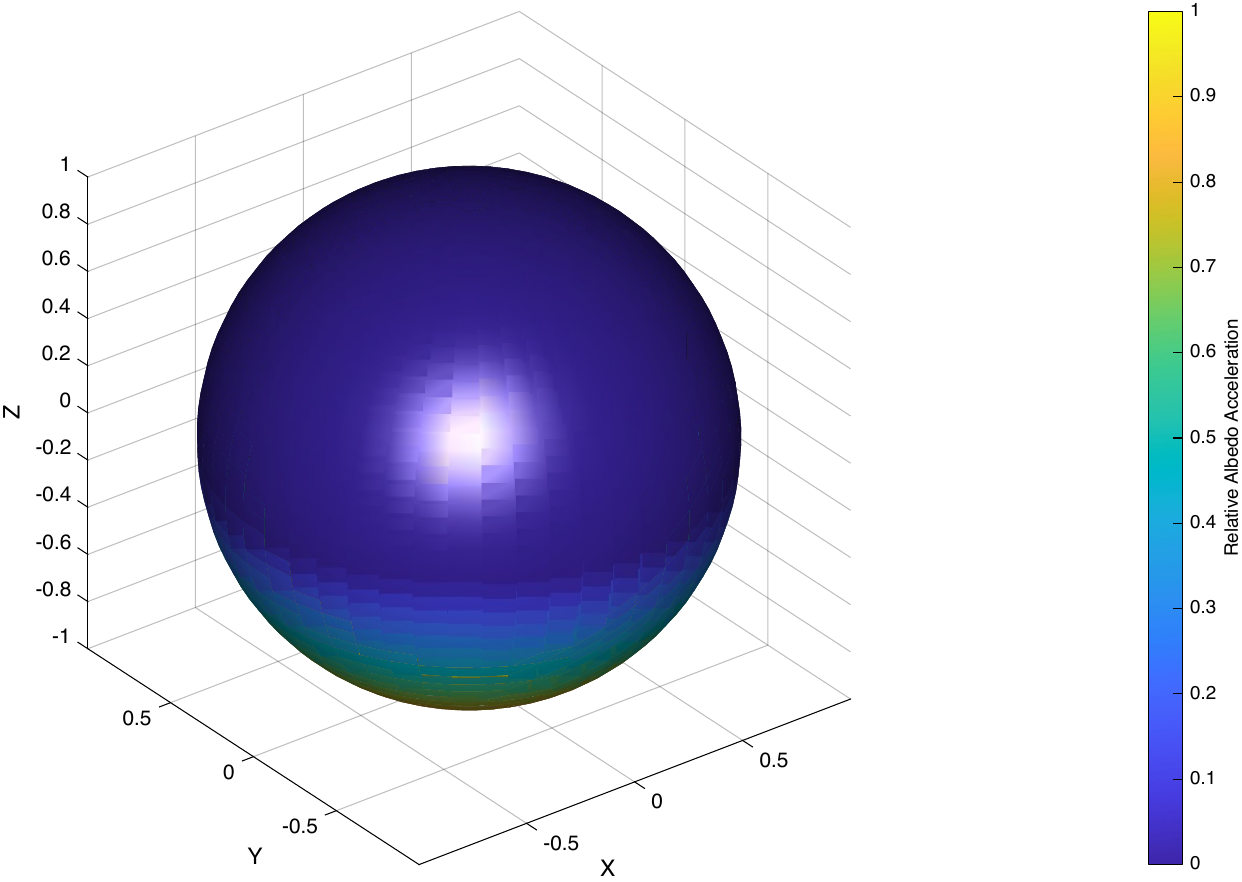}
\caption{Net albedo acceleration as a composite effect of varying pressures across all facets.}
\label{fig:sphere_eei}
\end{figure}

In a multi-facet spherical spacecraft model, the surface is divided into numerous small facets, each oriented differently relative to the Earth. Since Earth albedo is the sunlight reflected off the Earth’s surface, facets facing toward Earth receive higher reflected radiation pressure, resulting in greater acceleration, while those oriented away experience little to no albedo-induced force. This directional dependence, as depicted in Fig. \ref{fig:sphere_eei}, means the net acceleration on the spacecraft is a composite effect of varying pressures across all facets. Modeling the spacecraft as multiple facets allows precise calculation of localized albedo accelerations, improving accuracy in predicting its overall motion and attitude dynamics under the influence of Earth-reflected solar radiation. Ultimately, repeating this process for a formation of spacecraft will result in a higher-fidelity analysis of the Earth's Energy Imbalance.

%% file: Sections/5_conclusion.tex
\section{Conclusion}
\label{sec:conclusion}

This paper presents a method to achieve high-accuracy estimates of the Earth’s Energy Imbalance (EEI) through innovative use of distributed spacecraft. Utilizing high-fidelity modeling and precise control, a significant step forward is provided in addressing the limitations of existing methods. 
The central challenge of interpreting radiation flux-induced accelerations under external disturbance conditions necessitates a novel solution: the deployment of a multi-spacecraft formation mission to fill diurnal sampling gaps and achieve full global coverage. 
Key findings include a detailed comparison of two attitude control approaches, magnetotorquers and reaction wheels, in the context of this high-precision mission. Reaction wheels demonstrate superior accuracy for high-demand applications, while magneto-torquers proved effective for energy-limited scenarios, emphasizing the need for tailored control systems based on mission priorities. Future work intends to analyze relative perturbations at the facet level, and present detailed results on the EEI for science modeling.